\documentclass{amsart}

\usepackage{amsmath,amssymb,amsthm}
\newtheorem{theorem}{Theorem}
\newtheorem{lemma}{Lemma}

\author{Poo-Sung Park}
\address{Department of Mathematics Education\\ 
Kyungnam University\\
Changwon, 51767 \\
Republic of Korea}
\email{pspark@kyungnam.ac.kr}

\title[Multiplicative functions additive on sums of two squares]{Multiplicative functions which are\\ additive on sums of two nonzero squares}

\begin{document}

\thanks{This research was supported by Basic Science Research Program through the National Research Foundation of Korea(NRF) funded by the Ministry of  Science and ICT(NRF-2017R1A2B1010761).}

\begin{abstract}
Let $f$ be a multiplicative function which satisfies
\[
	f(a^2+b^2+c^2+d^2) = f(a^2+b^2)+f(c^2+d^2)
\]
for positive integers $a$, $b$, $c$, and $d$.
We show that $f$ is the identity function provided that $f(3)\,f(11) \ne 0$. Otherwise, $f(n)=0$ for all $n \ge 2$ except for $n=3,9,11$.
\end{abstract}

\maketitle

\section{Introduction}

Let $S$ be a set of arithmetic functions and $E$ be a set of positive integers. If an arithmetic function $f \in S$ is uniquely determined under the condition
\[
	f(m+n) = f(m)+f(n) \quad \text{for all $m,n \in E$},
\]
we call $E$ an \emph{additive uniqueness set for $S$}.

Claudia Spiro \cite{Spiro}, who coined the term, prove in 1992 that the set $P$ of all primes is an additive uniqueness set for the set of multiplicative functions $f$ with $f(p_0) \ne 0$ for some $p_0 \in P$. Since her igniting study, a number of mathematicians have been studying various problems related with the additive condition.

P. V. Chung \cite{Chung} in 1996 classified multiplicative functions satisfying
\[
	f(m^2+n^2) = f(m^2) + f(n^2) \quad \text{for all $m,n \in \mathbb{N}$}.
\]
He showed that the set of positive squares is an additive uniqueness set for completely multiplicative functions, but is not one for mere multiplicative functions.

He and B. M. Phong \cite{Chung-Phong} in 1999 proved that the set of triangular numbers is an additive uniqueness set for multiplicative functions. Also, so is the set of tetrahedral numbers.

In the present article, we consider the set of sums of two nonzero squares. That is, our question is whether a multiplicative function $f$ is the identity function or not provided $f$ satisfies
\[
	f(a^2+b^2+c^2+d^2) = f(a^2+b^2) + f(c^2+d^2) \quad \text{for all $a, b, c, d \in \mathbb{N}$}.
\]

Similar researches were performed by K.-H. Indlekofer and B. M. Phong \cite{Indlekofer-Phong} in 2006. They characterized all multiplicative function $f$ satisfying
\[
	f(n^2+m^2+k+1) = f(n^2+k) + f(m^2+1) 
\]
with fixed $k \in \mathbb{N}$ and all $n,m \in \mathbb{N}$.

Recently, in 2016, B. M. Phong \cite{Phong} studied multiplicative functions $f$ and $g$ satisfying
\[
	f(m^2+n^2+a+b) = g(m^2+a) + g(n^2+b) \quad \text{for all $m,n \in \mathbb{N}$},
\]
where $a$ and $b$ are some non-negative integers.

But, the condition of our problem can be said to be stronger, so the set of sums of two nonzero squares can be an additive uniqueness set for multiplicative functions $f$ with $f(n) \ne 0$ for some $n \in \mathbb{Z}\setminus\{1, 3, 9, 11\}$.

\section{Main theorem}

\begin{theorem}
If a multiplicative function $f$ satisfies
\[
	f(a^2+b^2+c^2+d^2) = f(a^2+b^2)+f(c^2+d^2)
\]
for positive integers $a$, $b$, $c$, and $d$, then $f$ is one of the following:
\begin{enumerate}
\item $f$ is the identity function,
\item $f(n) = 0$ for $n \ge 2$,
\item $f(3)\,f(9) \ne 0$ and $f(n) = 0$ for other $n \ge 2$,
\item $f(9) \ne 0$ and $f(n) = 0$ for other $n \ge 2$,
\item $f(11) \ne 0$ and $f(n) = 0$ for other $n \ge 2$.
\end{enumerate}

\end{theorem}

We need a condition for an integer to be represented as a sum of 4 nonzero squares to prove the main theorem. In 1911, E. Dubouis \cite{Dubouis} classified the general conditions.

\begin{lemma}[Dubouis]\label{Dubouis}
Every integer $n$ can be represented as a sum of $k$ nonzero squares except
\[
n = \begin{cases}
1, 3, 5, 9, 11, 17, 29, 41, 2\cdot4^m, 6\cdot4^m, 14\cdot4^m ~(m \ge 0) & \text{if } k = 4, \\
33 & \text{if } k = 5, \\
1, 2, \dots, k-1, k+1, k+2, k+4, k+5, k+7, k+10, k+13 & \text{if } k \ge 5.
\end{cases}
\]
\end{lemma}

Now, we compute some $f(n)$ for small $n$'s by using the following equalities:

\allowdisplaybreaks
\begin{align*}
f(1^2+1^2+1^2+1^2) 
	&= f(4) \\*
	&= 2\,f(2) \\
f(1^2+1^2+1^2+2^2)
	&= f(7) \\*
	&= f(2) + f(5) \\
f(1^2+1^2+2^2+2^2)
	&= f(2)\,f(5) \\*
	&= 2\,f(5) \\*
	&= f(2) + f(8) \\
f(1^2+1^2+1^2+3^2)
	&= f(3)\,f(4) \\*
	&= f(2) + f(2)f(5) \\
f(1^2+2^2+2^2+2^2)
	&= f(13) \\*
	&= f(5) + f(8) \\
f(1^2+1^2+2^2+3^2)
	&= f(3)\,f(5) \\*
	&= f(5) + f(2)f(5) \\
f(2^2+2^2+2^2+3^2)
	&= f(3)\,f(7) \\*
	&= f(8) + f(13)	
\end{align*}

Put $x = f(2)$, $y=f(3)$, and $z=f(5)$. Then, $f(4)=2x$. From the equalities for $f(10)$ we have that
\[
	xz = 2z = x+f(8).
\]

We divide two cases according to $z$. First, assume that $z \ne 0$. Then, $x=2$ and $f(8) = 2z-x$ from the above equality.

Since $f(4) = 2x$, from the equalities for $f(12)$ we obtain that
\[
	y \cdot 2x = x + xz \text{ or } 2y = 1+z.
\]
Now, from the equality $yz = z(1+x)$ for $f(15)$ we get $y = 3$ and $z = 5$.

Thus, if $f(5) \ne 0$, then we can easily check that $f(n) = n$ for $n = 1, 2, \dots, 21$. Note that $f(9)$, $f(11)$, and $f(17)$ can be calculated by using $f(2 \cdot 9)$, $f(2 \cdot 11)$, and $f(2 \cdot 17)$.

Let us consider the second case $z = 0$. Then, $f(8) = f(13) = -x$. From equalities for $f(12)$ and $f(21)$ we obtain that
\[
y \cdot 2x = x \text{ and }yx = -2x.
\]
Thus, $x$ should vanish and $f(n) = 0$ for $2 \le n \le 21$ and $n \ne 3,9,11,17,19$. But we can find that $f(17) = f(19)$ from $f(19) = f(17)+f(2)$.

In this case, we determine $f(17) = f(19) = f(25) = 0$ from
\[
	f(25) = f(8) + f(17) = f(5) + f(4)\,f(5).
\]
Also, we obtain $f(3)\,f(11) = 0$ from
\[
	f(33) = f(3)\,f(11) = f(8) + f(25) = f(13) + f(4)\,f(5).
\]

Now, we separate the proof into three cases.

\bigskip
\noindent (i) $f(n) = n$ for $1 \le n \le 21$

Now, we use induction to show that $f$ is the identity function in this case. Assume that $f(n) = n$ for all $n < N$. If $N = a^2+b^2+c^2+d^2$ for some positive integers $a$, $b$, $c$, and $d$, then $f(a^2+b^2) = a^2+b^2$ and $f(c^2+d^2) = c^2+d^2$ by induction hypothesis. Thus, $f(N) = N$.

If $N$ cannot be represented as a sum of four nonzero squares, $N$ is $29$, $41$, $2\cdot4^m$, $6\cdot4^m$, or $14\cdot4^m$ by Lemma \ref{Dubouis}.

We can obtain $f(29) = 29$ from 
\[
f(1^2+2^2+2^2+5^2) = f(2)\,f(17) = f(5)+f(29).
\]
Similarly, $f(41) = 41$ from 
\[
f(1^2+3^2+4^2+5^2) = f(3)\,f(17) = f(10)+f(41).
\]

If $N = 2\cdot4^m$, then $f(2\cdot4^m) = 4\,f(2\cdot4^{m-1}) = 2\cdot4^m$ from
\begin{align*}
f(4^{m-1}+4^{m-1}+4^m+4^m) 
	&= f(2\cdot4^{m-1})\,f(5) \\
	&= f(2\cdot4^{m-1}) + f(2\cdot4^m)
\end{align*}
and induction hypothesis.

Other cases can be calculated by
\[
f(6 \cdot 4^m) = f(3)\,f(2\cdot4^m)
\text{ and }
f(14\cdot4^m) = f(7)\,f(2\cdot4^m).
\]

Hence, $f(n) = n$ for all $n \ge 1$.

\bigskip
\noindent (ii) $f(3) = 0$

We have $f(n) = 0$ for $2 \le n \le 21$ and $n \ne 9, 11$. We use induction again. Assume that $f(n)=0$ for $n < N$ except for $n = 9,11$.

Suppose $N$ is a sum of four nonzero squares. Note neither $9$ nor $11$ is a sum of two nonzero squares. So $f(9)$ and $f(11)$ are not used to calculate $f(N)$. Thus, $f(N)=0$. Especially, 
\[
	f(99) = f(9)\,f(11) = f(1^2+1^2) + f(4^2+9^2) = 0
\]
and thus, $f(9)=0$ or $f(11)=0$.

If $N$ is not a sum of four nonzero squares, then $N$ is $29$, $41$, $2\cdot4^m$, $6\cdot4^m$, or $14\cdot4^m$.  We can obtain $f(N)$ for these $N$ by the same ways as the previous case.

\bigskip
\noindent (iii) $f(11) = 0$

In this case, we have that $f(n)=0$ for $2 \le n \le 21$ and $n \ne 3, 9$. By induction, $f(N)=0$ if $N$ is a sum of four nonzero squares. The exceptional numbers $29$, $41$, $2\cdot4^m$, $6\cdot4^m$, and $14\cdot4^m$ can be dealt with the same ways as the previous case. But, since $3$ and $9$ are not sums of two nonzero squares and they are not relatively prime, we cannot determine $f(3)$ and $f(9)$.

\end{document}